\documentstyle[12pt]{amsart}

\makeatletter                                            
\renewcommand{\@begintheorem}[2]{                        
\rm \trivlist \item [\hskip \labelsep {\bf #2\ \ #1.}]   
                                }                        
\makeatother                                             

\newcommand{\newsubsection}%
{{\noindent\bf\refstepcounter{subsection}\thesubsection\ \ }}

\newcommand{\newsubsubsection}%
{{\bf\refstepcounter{subsubsection}\thesubsubsection\ \ }}

\newcommand{\ts}{\vspace{\baselineskip}\noindent{\bf Proof.}$\;\;$}
\newcommand{\ZZ}{{\bf Z}}
\newcommand{\QQ}{{\bf Q}}

\newcommand{\CC}{{\bf C}}

\newcommand{\FF}{{\bf F}}

\newcommand{\PP}{{\bf P}}

\newcommand{\cC}{{\cal C}}

\newcommand{\cL}{{\cal L}}
\newcommand{\cM}{{\cal M}}

\newcommand{\cO}{{\cal O}}

\newcommand{\cT}{{\cal T}}
\newcommand{\cX}{{\cal X}}

\def\DynkinEEE#1#2#3#4#5#6#7
   {\begin{picture}(124, 50)%
   \put(5, 35){\circle{4}}%
   \put(7, 35){\line(1, 0){28}}%
   \put(37, 35){\circle{4}}%
   \put(39, 35){\line(1, 0){28}}%
   \put(69, 35){\circle{4}}%
   \put(71, 35){\line(1, 0){28}}%
   \put(101, 35){\circle{4}}%
   \put(103, 35){\line(1, 0){28}}%
   \put(133, 35){\circle{4}}%
   \put(69, 33){\line(0, -11){28}}%
   \put(69, 3){\circle{4}}%
   \put(5, 41){\makebox(0, 0)[b]{$#1$}}%
   \put(37, 41){\makebox(0, 0)[b]{$#3$}}%
   \put(69, 41){\makebox(0, 0)[b]{$#4$}}%
   \put(101, 41){\makebox(0, 0)[b]{$#5$}}%
   \put(133, 41){\makebox(0, 0)[b]{$#6$}}%
   \put(69, 0){\makebox(0, 0)[t]{$#2$}}%
   \put(-60, 35){\makebox(0, 0)[lb]{$#7$}}%
\end{picture}}

\begin{document}

\title{A family of marked cubic surfaces and \\
the root system $D_5$}

\author{Elisabetta Colombo}
\address{Dipartimento di Matematica, Universit\`a di Milano,
  via Saldini 50, I-20133 Milano, Italia}
\author{Bert van Geemen}

\begin{abstract}
We define and study a family of cubic surfaces in the projectivized
tangent bundle over a four dimensional projective space associated to the
root system $D_5$.
The $27$ lines are rational over the base and we determine the classifying map
to the moduli space of marked cubic surfaces. This map has degree two
and we use it to get short proofs for some results on the Chow group
of the moduli space of marked cubic surfaces.
\end{abstract}

\maketitle

A marked cubic surface is a smooth cubic surface $S$ with an ordered set of six
skew lines on $S$.
The symmetry group of the configuration of the $27$ lines on a smooth
cubic surface is the Weyl group of the root system $E_6$. This group
acts on the moduli space of marked cubic surfaces $\cM$.
In various constructions of $\cM$
certain subgroups of $W(E_6)$ play a prominent role.

Starting with the description of a marked cubic surface as the blow-up of
the projective plane in six ordered points the action of the symmetric group
$S_6\subset W(E_6)$ on $\cM$ is very evident.
In this approach, $\cM$ is constructed as a double cover of a projective
space branched along a quartic threefold, cf.\ \cite{DGK}, 2.11.
The universal family of marked cubic surfaces is given there by the
Cremona hexaedral form of a cubic surface (cf.\ \cite{Co1}, section 4,
\cite{Dcr}, section 3), this family has an obvious $S_6$-action.

In his construction of $\cM$
and its desingularisation $\cC$, Naruki \cite{Naruki} used the subgroup
$W(F_4)$, which is associated to a tritangent plane of a cubic surface.
In his approach, $\cM$ is a compactification of the maximal torus of a complex
Lie group of type $D_4$. The universal family is given explicitly, but it is
rather complicated.

In this paper we consider the subgroup $W(D_5)$ of $W(E_6)$.
It corresponds to the stabiliser of a line on a cubic surface.
The standard five dimensional represention of $W(D_5)$ gives an action of
this Weyl group on projective four space. We define a family of marked cubic
surfaces $\cX$ in the projectivized tangent bundle of this projective space.
The family is defined by a very simple determinantal equation
and $W(D_5)$ acts on it.
The 27 lines and the 45 tritangent planes are also given by simple
expressions.

Naruki's construction of $\cM$ and subsequent work of Allcock and Freitag on
an explicit projective embedding of $\cM$ into $\PP^9$
(using modular forms on a four ball),
allows us to find a vector space of quintic polynomials
in five variables which gives the classifying map from an open subset
of the projective four space to
$\cM\subset \PP^9$. We find that the classifying map has degree two
and that the covering involution is simply inversion of the (natural) coordinates.
We extend the classifying map
to an explicit blow-up of the projective space. This allows us to recover some
results on the Chow group of the moduli space of marked cubic surfaces which
we obtained earlier in \cite{CvG} using Naruki's model.

\section{The root system $D_5$}

\subsection{} We recall some basic facts on root systems and cubic surfaces.
For root systems we refer to \cite{Hu}, for cubic surfaces to \cite{DGK}
and the references given there.

\subsection{The root lattice $Q(D_5)$}
Let $Q(D_5)\cong\ZZ^5$ be the root lattice of type $D_5$.
It can be realised in the complex vector space of linear
forms on $\CC^5$ as
$$
Q(D_5)=\{\sum_{i=1}^5 a_ix_i:\;
  a_i\in\ZZ,\;a_1+\ldots+a_5\equiv 0\;{\rm mod}\,2\},
\qquad x_i\cdot x_j=\delta_{ij}
$$
where `$\cdot$' indicates the bilinear form on
$Q_\CC:=Q(D_5)\otimes_\ZZ\CC$. A $\ZZ$-basis of $Q(D_5)$ of simple roots
is given by:
$$
\alpha_1=h_{12}=x_1-x_2,\quad\alpha_2=h_{123}=x_4+x_5,\quad
\alpha_3=h_{23}=x_2-x_3,\quad\ldots\quad \alpha_5=h_{45}=x_4-x_5.
$$
The roots of $D_5$ are the $\pm x_i\pm x_j\in Q(D_5)$, $i\neq j$,
the positive roots are $h_{ij}=x_i-x_j$, $1\leq i<j\leq 5$,
and $h_{klm}=x_i+x_j$ where
$\{i,\ldots,m\}=\{1,\ldots,5\}$.

\subsection{The Weyl group $W(D_5)$}\label{wd5}
The Weyl group $W(D_5)$ is the subgroup of $GL(Q_\CC)$ generated
by the reflections in the hyperplanes perpendicular to the
roots of $D_5$. The reflection defined
by $h_{ij}$ is denoted by $s_{ij}$. It permutes $x_i$ and $x_j$ and fixes the
other $x_k$, this gives a subgroup $S_5\subset W(D_5)$.
The reflection $s_{klm}$ defined by $h_{klm}$ maps $x_i$
to $-x_j$ and $x_j$ to $-x_i$ and fixes $x_k,x_l,x_m$.

\subsection{The root system $E_6$ and cubic surfaces}\label{lines}
The Picard group of a smooth complex cubic surface $S$ is a free
$\ZZ$-module of rank $7$ which has a basis $a_0,a_1,\ldots a_6$
such that the
intersection form is given by $x_0^2-(x_1^2+\ldots+x_6^2)$.
The canonical class is $k=-3a_0+(a_1+\ldots+a_6)$ and $k^\perp$ is isomorphic,
up to sign,
to the root lattice $Q(E_6)$. The classes in $Pic(S)$
of the $27$ lines on $S$ are the elements $l\in Pic(S)$ with
$l^2=-1$ and $k\cdot l=-1$. In particular, they are not in $Q(E_6)$,
but they project onto a $W(E_6)$-orbit of $27$ weights in $Q(E_6)\otimes\QQ$.
The (classes of the) lines are
$$
a_i,\qquad b_j=2a_0-(a_1+\ldots+\widehat{a_j}+\ldots +a_6),\qquad
c_{kl}=a_0-(a_k+a_l),
$$
with $i,j,k,l\in\{1,\ldots,6\},\;k\neq l$.
The six lines $a_1,\ldots,a_6$ are skew and can be blown down
to give a birational  morphism $S\rightarrow \PP^2$ which maps
these lines to six points $p_i\in\PP^2$. The line  $b_i$
maps to the conic through the points
$p_1,\ldots,\widehat{p_i},\ldots,p_6$, $1\leq i\leq 6$, and the
line $c_{ij}=c_{ji}$ maps to the line
through $p_i$ and $p_j$.
The Weyl group of $E_6$ is generated by reflections defined by
the $-2$-classes in $Q(E_6)\subset Pic(S)$, in particular,
$W(E_6)$ acts by isometries on $Pic(S)$.

\subsection{Orbits of $W(D_5)$ on the 27 lines}\label{orbits lines}
We identify $W(D_5)\subset W(E_6)$ with the stabiliser of the line $a_6$.
Equivalently, $D_5\subset E_6$ is generated by the simple roots
$\alpha_1,\ldots ,\alpha_5$ in the $E_6$-diagram below.
$$
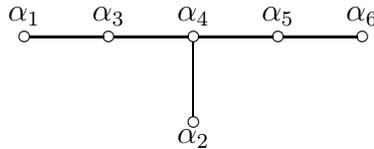
\DynkinEEE{\alpha_1}{\alpha_2}{\alpha_3}{\alpha_4}{\alpha_5}{\alpha_6}
$$

The subgroup  $S_5\subset W(D_5)$ acts on the lines by
permutation of the indices. The reflection $s_{123}$,
whose action on the six points is given by the
Cremona transformation in $p_1,p_2,p_3$, maps
$$
s_{123}:\qquad a_1\leftrightarrow c_{23},
\quad a_2\leftrightarrow c_{13},\quad
a_3\leftrightarrow c_{12},\quad
b_4\leftrightarrow c_{56},\quad b_5\leftrightarrow c_{46},\quad
b_6\leftrightarrow c_{45}
$$
and fixes the other lines.
Thus the $W(D_5)$-orbits in the set of lines are:
$$
\{a_6\},\qquad \{a_1,\ldots,a_5,c_{12},\ldots,c_{45},b_6\},\qquad
\{b_1,\ldots,b_5,c_{16},\ldots,c_{56}\}.
$$
Note that the 5 pairs of lines $(b_1,c_{16}),\ldots,(b_5,c_{56})$
meet $a_6$ and the other $27-1-10=16$ lines are disjoint from $a_6$.

\subsection{Orbits of $W(D_5)$ on the 45 tritangent planes}
\label{orbits tritangents}
A tritangent of a smooth cubic surface $S$ is a plane section of
$S$ which consists of three lines. As the embedding $S\hookrightarrow \PP^3$
is anticanonical, the tritangents correspond to sets of three lines whose
sum, in $Pic(S)$, is $-k$.

There are 45 tritangent planes,
they correspond to the following sets of lines:
$$
\{a_i,b_j,c_{ij}\},\qquad \{c_{ij},c_{kl},c_{mn}\}
$$
where $i\neq j$ and $\{i,j,\ldots,n\}=\{1,2,\ldots,6\}$.
There are two $W(D_5)$-orbits on the set of tritangent divisors,
of length 5 and 40 respectively, they are:
$$
\{\ldots,\{a_6,b_i,c_{i6}\},\ldots\}_{1\leq i\leq 5},\qquad
\{\ldots,\{a_p,b_q,c_{pq}\},\ldots,\{c_{ij},c_{kl},c_{mn}\},\ldots\},
$$
where in the last set we have $p\neq q$, $1\leq p\leq 5$, $1\leq q\leq 6$ and
$\{i,j,\ldots,n\}=\{1,2,\ldots,6\}$.

\section{A family of marked cubic surfaces}

\subsection{} We will construct a family of marked cubic surfaces
embedded in the projectivised
tangent bundle over a $\PP^4$ naturally associated to the
root system $D_5$.

Recall that $Q_\CC:=Q(D_5)\otimes\CC\cong\CC^5$. Let ${\cal T}$
be the tangent bundle of the projective space $\PP(Q_\CC)\cong\PP^4$.
The fibre of ${\cal T}\rightarrow \PP(Q_\CC)$
over $<x>\in \PP(Q_\CC)$ is the vector space
$Q_\CC/\langle x\rangle$.

Let $\PP({\cal T})\rightarrow \PP(Q_\CC)$ be
the associated projective space bundle,
its fibres are $\PP^3$'s.
The Euler sequence, which maps the trivial bundle $\PP(Q_\CC)\times Q_\CC$
over $\PP(Q_\CC)$ to ${\cal T}(-1)$, gives a natural map $\pi$:
$$
\begin{array}{ccc}
\PP(Q_\CC)\times (Q_\CC-\{0\})&\stackrel{\pi}{\longrightarrow}& \PP({\cal T})\\
\downarrow &&\downarrow\\
\PP(Q_\CC)&=&\PP(Q_\CC).
\end{array}
$$
We will denote the homogeneous coordinates on $\PP(Q_\CC)$ by $x_i$ and
the $X_i$ are the coordinates on $Q_\CC$.
A family of surfaces ${\cal X}\subset \PP({\cal T})$ over $\PP(Q_\CC)$
is determined by the polynomial $F$ in the $x_i$ and $X_i$
which defines the inverse image of
${\cal X}$ in $\PP(Q_\CC)\times (Q_\CC-\{0\})
\cong\PP^4\times(\CC^5-\{0\})$.

\subsection{The determinant}\label{definition F}
Consider the following polynomial in $\CC[x_1,\ldots,x_5,X_1,\ldots,X_5]$,
bihomogeneous of degree $(7,3)$ which is defined as the determinant of the
following matrix:
$$
F(x,X)=\det\left(
\begin{array}{rrrrr}
1&1&1&1&1\\
x_1^2&x_2^2&x_3^2&x_4^2&x_5^2\\
x_1^4&x_2^4&x_3^4&x_4^4&x_5^4\\
x_1X_1&x_2X_2&x_3X_3&x_4X_4&x_5X_5\\
X_1^2&X_2^2&X_3^2&X_4^2&X_5^2\\
\end{array}\right).
$$
Thus we have:
$$
F(x,X)=\sum_{i\neq j} f_{ij}(x)X_iX_j^2,\qquad
f_{12}=-x_1(x_3^2-x_4^2)(x_3^2-x_5^2)(x_4^2-x_5^2),
$$
and the other polynomials $f_{ij}$ are determined by:
$$
f_{\sigma^{-1}(i)\sigma^{-1}(j)}(x_1,\ldots,x_5)=
sign(\sigma)f_{ij}(x_{\sigma(1)},\ldots,x_{\sigma(5)}),
$$
where $\sigma\in S_5$.

\subsection{Descend to the tangent bundle}\label{definition U}
We define an open subset of $\PP(Q_\CC)$ by:
$$
U:=\,
\{(x_1:\ldots:x_5)\in\PP(Q_\CC): \;
x_i\neq \pm x_j \text{ for}\; i\neq j\,\}.
$$
The following lemma implies that for $x\in U$ the polynomial $F$ from section
\ref{definition F} defines
the cone over a cubic surface in $(Q_\CC-\{0\})/\CC^*\cong\PP^4$ with vertex
$<x>$.

\subsection{Lemma}\label{descend}
There is a family of cubic surfaces
${\cal X}\hookrightarrow \PP({\cal T})_{|U}$
over $U$ whose inverse image in $U\times (Q_\CC-\{0\})$
under $\pi$ is defined by $F(x,X)=0$.

\ts
Note that for $x\in U$ the polynomial $F(x,X)$
in the $X_i$ is not identically zero.
We verify that $F(x,X+tx)=F(x,X)$ for all $t$,
so we get a well defined family of cubic
surfaces ${\cal X}$
in $\PP({\cal T}_Q)$, i.e.\ for $x\in \PP(Q_\CC)$ the surface
${\cal X}_x\subset\PP(Q_\CC/<x>)$ is
defined by $F(x,X)=0$.

Substituting $X_i:=X_i+tx_i$, $1\leq i\leq 5$, in the matrix defining
$F$, one easily sees that the last two rows become
linear combinations of the rows of the matrix and that the determinant
doesn't change.
\qed

\subsection{The $W(D_5)$-action on the family}
The subset $U$ is invariant under the action of $W(D_5)$ and this action
lifts to an action of $W(D_5)$ on ${\cal T}_{|U}$.
In fact, $g\in W(D_5)$ acts on the $x_i$
via the standard representation. As this action is linear, the action
on ${\cal T}$ is the one induced by the standard representation on the $X_i$.
For $g\in W(D_5)$ and $(x,X)\in \PP(Q_\CC)\times (Q_\CC-\{0\})$ we thus have
$g(x,X):=(gx,gX)$.

The group $W(D_5)$ is generated by $S_5$ and
$s_{123}$, which maps $x_4\mapsto -x_5$, $x_5\mapsto -x_4$
and fixes the other $x_i$.
It is then easy to check from the determinant defining $F$ that
$F(gx,gX)=\det(g)F(x,X)$. Thus
the family ${\cal X}$ is invariant under the action of $W(D_5)$ on
$\PP({\cal T})_{|U}$.

\subsection{Determinants}\label{determinant}
To find lines in these cubic surfaces we need the following
determinantal identity:
$$
\det\left(
\begin{array}{ccccc}
1&\ldots&1&\ldots&1\\
x_1&\ldots&x_i&\ldots&x_n\\
\vdots&\ldots&\vdots&\ldots&\vdots\\
x_1^{n-2}&\ldots&x_i^{n-2}&\ldots&x_n^{n-2}\\
x_2\cdots x_{n-1}&\ldots&x_1\cdots\widehat{x_i}\cdots x_n
&\ldots&x_1\cdots x_{n-1}\\
\end{array}\right)
=(-1)^{(n-1)(n-2)/2}\prod_{i<j} (x_i-x_j).
$$
This holds because the determinant
is a homogeneous polynomial of degree $n(n-1)/2$ which is divisible by
$x_i-x_j$ for every $i\neq j$ and the coefficients of
$x_1x_2^2\ldots x_{n-1}^{n-1}$ on both sides are equal.
Note that the
Van der Monde determinant is obtained by replacing the last row
by $x_1^{n-1},\ldots,x_n^{n-1}$
and that it equals $(-1)^{n(n-1)}\prod_{i<j} (x_i-x_j)$.
In particular, the sum of the Van der Monde determinant and the one
above is zero if $n=4$.

\subsection{Lines in ${\cX}$}
A line in $\PP({\cal T}_x)$ is defined by a two dimensional subspace
of ${\cal T}_x=Q_\CC/\langle x\rangle$ and hence
by a three dimensional subspace of $Q_\CC$ which
contains $x$. In the following theorem we give some lines on ${\cal X}_x$
by specifying these three dimensional subspaces as $\langle x,s(x),t(x)\rangle$
with $s(x),t(x)\in Q_\CC\cong \CC^5$.

The Weyl group $W(D_5)$ acts on ${\cal X}$.
As $g\in W(D_5)$ maps a section $\sigma$ with $\sigma(x)=(x,s(x))$ to
$(gx,gs(x))$,
we get an action of $W(D_5)$ on the sections
with $g(\sigma)(x)=(x,gs(g^{-1}x))$.

Recall that a double six on a cubic surface is a set of twelve lines
$\{a_1,\ldots,a_6,b_1,\ldots,b_6\}$ such that
the intersection products are
(cf.\ \cite{Hunt}):
$$
a_i\cdot a_j=0,\quad b_i\cdot b_j=0,\quad a_i\cdot b_j=1, \quad a_i\cdot b_i=0,
\qquad \text{for}\;i\neq j.
$$

\subsection{Theorem}\label{double six}
For $x\in U$, the cubic surface ${\cal X}_x$ has a double six formed by
the following twelve lines:
$$
a_1=\langle x,
(-1,1,1,1,1),\, (-x_1^2,x_2^2,x_3^2,x_4^2,x_5^2)\,\rangle,
$$
the lines $a_2,\ldots,a_5$ defined as $a_i=s_{1i}(a_1)$ with
$s_{1i}\in S_5\subset W(D_5)$,
$$
a_6=\langle x,\,
(x_1^{-1},x_2^{-1},x_3^{-1},x_4^{-1},x_5^{-1}),\,
(x_1^3,x_2^3,x_3^3,x_4^3,x_5^3)\,\rangle,
$$
$$
b_1=\langle x,\,
(1,0,0,0,0),
(0,x_2^3-x_3x_4x_5,x_3^3-x_2x_4x_5,x_4^3-x_2x_3x_5,x_5^3-x_2x_3x_4)
\,\rangle,
$$
the lines $b_2,\ldots,b_5$ defined as $b_i=s_{1i}(b_1)$ with
$s_{1i}\in S_5\subset W(D_5)$, and
$$
b_6=\langle x,
(1,1,1,1,1),\, (x_1^2,x_2^2,x_3^2,x_4^2,x_5^2)\,\rangle.
$$

\ts
In view of Lemma \ref{descend}, to check that $\langle x,s(x),t(x)\rangle$
defines a line in $F(x,X)=0$ it suffices to check that
$F(x,s(x)+\lambda t(x))=0$ for all $\lambda$.

Thus $b_6$ lies in $F(x,X)=0$ since
upon substituting $X_i:=1+\lambda x_i^2$ in the
matrix defining $F$, the last row becomes a linear combination
of the first three hence its determinant is zero.

Recall that $g\in W(D_5)$ acts on sections by
$g:(x,s(x)) \mapsto (x,gs(g^{-1}x))$ and maps lines to lines.
Thus $b_6$ is fixed under the $s_{ij}$
(which permute $x_i$ and $x_j$ and fix the other coordinates).
Its image under $s_{123}$ (which maps $x_4\mapsto -x_5$,
$x_5\mapsto -x_4$ and fixes the other $x_i$) is:
$$
c_{45}=s_{123}(b_6)=\langle x,
(1,1,1,-1,-1),\, (x_1^2,x_2^2,x_3^2,-x_4^2,-x_5^2)\,\rangle.
$$
Applying $s_{145}$ (which maps $x_2\mapsto -x_3$,
$x_3\mapsto -x_2$ and fixes the other $x_i$) we get:
$$
a_1=s_{145}(c_{45})=\langle x,
(-1,1,1,1,1),\, (-x_1^2,x_2^2,x_3^2,x_4^2,x_5^2)\,\rangle.
$$
As the $a_i$, $1\leq i\leq 5$, are in the $W(D_5)$-orbit of $b_6$, these
lines are also on the cubic surface ${\cal X}_x$.

An argument similar to the one above shows that $a_6$
is a line on ${\cal X}_x$.
Under the $W(D_5)$-action this line is invariant. In fact,
all coordinate functions
are of odd degree so sign changes are merely reparametrisations of
the same line.

The linearity of the determinant in each row
of a matrix and the identity given in \ref{determinant}
easily show that $b_1$ lies on ${\cal X}_x$, and so the $b_i$, $1\leq i\leq 5$
are on ${\cal X}_x$. Thus the twelve lines in the theorem are on ${\cal X}_x$.

Two lines
$l=\langle x,v_1,v_2\rangle,m=\langle x,w_1,w_2\rangle$ in $\PP({\cal T}_x)$
are skew iff the $5\times 5$ matrix
whose rows are $x,v_1,v_2,w_1,w_2$ has maximal rank.
In particular,
the lines $b_6$ and $a_6$ are skew because the $5\times 5$ matrix
they define is, up to multiplication of the $i$-th
column by $x_i$, a Van der Monde matrix. Using the $W(D_5)$-action on the
lines (this action fixes $a_6$) we find that
each of the 16 lines in the $W(D_5)$-orbit of $b_6$ is skew with $a_6$.
In particular, $a_6$ is skew with the lines denoted by $a_1,\ldots,a_5$ above.
It is also easy to see that $a_i$ and $a_j$ are skew if $1\leq i< j\leq 5$.
Hence the six $a_i$ are skew.

The lines $b_i$ and $b_j$ are skew if $1\leq i< j\leq 5$, in fact in case
$i=1,j=2$ the
$5\times 5$ matrix now has maximal rank since:
$$
\det\left(
\begin{array}{ccc} x_3&x_4&x_5\\
x_3^3-x_2x_4x_5 & x_4^3-x_2x_3x_5 & x_5^3-x_2x_3x_4\\
x_3^3-x_1x_4x_5 & x_4^3-x_1x_3x_5 & x_5^3-x_1x_3x_4
\end{array}\right)=(x_1-x_2)(x_3^2-x_4^2)(x_3^2-x_5^2)(x_4^2-x_5^2),
$$
and the other cases are similar. The lines $b_i,b_6$, $1\leq i\leq 5$, are
skew since the $5\times 5$ matrix (essentially a $4\times 4$ matrix)
can be seen to have maximal rank using the remarks in section \ref{determinant}.
Hence the six $b_i$ are skew.

We already observed that $b_6$ and $a_6$ are skew. It is not hard
to verify that
$a_i$ and $b_i$ are skew.
We show that $a_i$ and $b_j$ intersect if $i\neq j$.
The line $b_6$ meets $a_1$ in the point:
$$
\begin{array}{rcl}
(0,x_1^2-x_2^2,\ldots,x_1^2-x_5^2)
&=&x_1^2(1,1,1,1,1)-(x_1^2,x_2^2,x_3^2,x_4^2,x_5^2)\\
&=&x_1^2(-1,1,1,1,1)-(-x_1^2,x_2^2,x_3^2,x_4^2,x_5^2).
\end{array}
$$
Similarly, $b_6$ meets $a_i$ for $i=2,\ldots,5$.
The lines $a_1$ and $b_5$ meet, since the $5\times 5$ matrix reduces to a
$4\times 4$ matrix which, upon substituting $x_1:=-x_1$ and multiplying
the first column by $-1$, is of the type discussed in section \ref{determinant}.
Applying  suitable elements $W(D_5)$ we get that $a_i$ and $b_j$ meet
if $i\neq j$ and $1\leq i,j\leq 5$.
Finally one verifies with similar arguments that $a_6$ meets $b_i$ if
$1\leq i\leq 5$.
\qed

\subsection{} Using the information on the lines we found in
the general ${\cal X}_x$, it is not hard to show that this surface must be smooth.
All we will need is that for $x\in U$, the cubic surface ${\cal X}_x$
has a double six.

\subsection{Lemma} The cubic surface ${\cal X}_x$ is smooth if $x\in U$, that is,
$x=(x_1:\ldots:x_5)$ and $x_i\neq \pm x_j$ for $i\neq j$.

\ts
Let $x\in U$, then ${\cal X}_x$ has six skew lines, the $a_i$,
and for each $j\in \{1,\ldots, 6\}$
there is a line, $b_j$, which meets all $a_i$ except $a_j$.
This easily implies that ${\cal X}_x$ is irreducible.
Suppose that $p\in {\cal X}_x$ is a singular point.
We choose coordinates such that $p=(0:0:0:1)$. Then
the equation defining ${\cal X}_x$ is
$G(x,y,z)t+H(x,y,z)=0$ for some $G,H\in\CC[x,y,z]$
homogeneous of degree two and three respectively.

If $G=0$, the surface ${\cal X}_x$ is a cone over the plane curve
defined by $H=0$, which must be irreducible,
thus all lines on ${\cal X}_x$ pass through $p$, a contradiction.
If $H=0$ the surface would be reducible, so we conclude that that both $G$
and $H$ are non-zero and have no common factor.
Projection from $p$ to the plane $t=0$ gives a birational isomorphism
${\cal X}_x\rightarrow \PP^2$. The inverse is given by $q=(x:y:z)\mapsto
(G(q)x:G(q)y:G(q)z:-H(q))$, which is an isomorphism on the complement
of the finite set $B=(G=0)\cap(H=0)$.
The lines containing $p$ project to the points in $B$, hence there are at most
six such lines.
Let $l$ be a line on
${\cal X}_x$ not containing $p$. The plane spanned by $p$ and $l$
intersects ${\cal X}_x$ in $l$ and a conic containing the singular point $p$.
Hence this conic is reducible and consists of two lines
(possibly a double line) passing through $p$. If there are no double lines,
the lines on ${\cal X}_x$ project to the points of $B$ or to the lines
spanned by two points in $B$. It is now easy to see that ${\cal X}_x$ does not
contain six skew lines.

Finally assume there is a plane $P$
such that there are only two lines $l,m$ with
$P\cap {\cal X}_x=l\cup m$, possibly $l=m$.
As each line in ${\cal X}_x$ intersects $P$ (or lies in $P$), there are at
least three of the six $a_i$'s which properly
intersect one of the two lines $l,m$.
Call that line $l$ and let $a_1,a_2,a_3$ be three $a_i$'s meeting $l$
in points. The lines $b_4,b_5,b_6$ all meet
$a_1,a_2,a_3$.
In case $b_j=l$ for one such $b_j$, say $j=6$,
the lines $a_6,b_1,\ldots,b_5$ must intersect $m$ properly.
The plane spanned by $a_6$ and $m$ then meets at least four of the five lines
$b_1,\ldots,b_5$ each in at least
two distinct points, so these lines must lie in that plane and hence they meet,
a contradiction, so $b_j\neq l$ for $j=4,5,6$.
Three skew lines in $\PP^3$ are contained in a unique smooth quadric,
let $Q$ be the quadric containing $a_1,a_2,a_3$.
Then $b_4,b_5,b_6$ also lie in $Q$ since each meets any of these three $a_i$
and similarly $l\subset Q$.
Thus $Q\cap{\cal X}_x$ contains at least $7$ distinct lines, but
an irreducible cubic surface
can meet a quadric in at most $6$ lines.
We conclude that a cubic surface having a double six must be smooth.
\qed

\subsection{The 27 lines}\label{27 lines}
A smooth cubic surface contains 27 lines. For $x\in U$, we already found a double
six $a_1,\ldots,b_6$ in ${\cal X}_x$.
The remaining 15 lines are the $c_{ij}=c_{ji}$,
$1\leq i<j\leq 6$. The line $c_{ij}$ is the third line in the intersection of
the plane spanned by $a_i$ and $b_j$ with ${\cal X}_x$. Alternatively,
one can find
the $c_{ij}$ using the $W(D_5)$ action on the lines in ${\cal X}_x$.
For example, $c_{45}$ was given in the proof of Theorem \ref{double six} and
$$
c_{16}=\langle x,
(1,0,0,0,0),\,
(0,x_2^3+x_3x_4x_5,x_3^3+x_2x_4x_5,x_4^3+x_2x_3x_5,x_5^3+x_2x_3x_4)\,\rangle.
$$

\section{Tritangent planes and cross ratios}
\subsection{Tritangent planes}\label{tritangents}
The lines $a_6,b_1, c_{16}$ lie on the tritangent plane:
$$
t_{61}:\qquad\sum_{j=1}^5 b_{1j}X_j:=
\det\left(
\begin{array}{rrrr}
 1&1&1&1\\
x_2^2&x_3^2&x_4^2&x_5^2\\
x_2^4&x_3^4&x_4^4&x_5^4\\
x_2X_2&x_3X_3&x_4X_4&x_5X_5\\
\end{array}\right)=0,
$$
Explicitly, we get
$$
b_{11}=0,\quad b_{12}=x_2(x_3^2 - x_4^2)(x_3^2 - x_5^2)(x_4^2 - x_5^2),\ldots
$$

The lines $a_1,b_6, c_{16}$ lie on the tritangent plane:
$$
t_{16}:\qquad\sum_{j=1}^5 a_{1j}X_j:=
\det\left(
\begin{array}{rrrr}
 1&1&1&1\\
x_2&x_3&x_4&x_5\\
x_2^2&x_3^2&x_4^2&x_5^2\\
X_2&X_3&X_4&X_5\\
\end{array}\right)=0,
$$
Explicitly, we get
$$
a_{11}=0,\quad a_{12}=(x_3 - x_4)(x_3 - x_5)(x_4 - x_5),
\ldots .
$$
Since the determinant is linear in each row,
a line $l=\langle x,v,w\rangle$ lies in $t_{ij}$ if upon substituting
$X:=(X_1,\ldots,X_5)=x,v,w$ the determinant defining $t_{ij}$ becomes zero.
To check that $c_{16}$ lies on $t_{16}$ one can use the linearity of the
matrix in the last row and the identity \ref{determinant}
to see that the two determinants cancel.

The set  of tritangent planes is the union of two $W(D_5)$-orbits. The
tritangent planes $t_{16},t_{61}$ are in
distinct orbits, so we essentially determined all tritangent planes.

\subsection{Cross ratios}\label{cross ratios}
Given a tritangent plane $t$ of a smooth cubic surface $X$ and one of the
three lines $l\subset t\cap X$, there are 4 other tritangent
planes $t_1,\ldots,t_4$ of $X$ which contain $l$. These four planes give
four points in the $\PP^1$ of planes containing $l$. Any one of the six
cross ratios of these 4 points is called a cross ratio associated to $X$ and $t$.
In fact, Cayley already observed that
these six cross ratios do not depend on the choice of the line
$l\subset t\cap X$, cf.\ \cite{Naruki}.

The pencil of lines on $l$ is spanned by $t_1$ and $t_2$
and if we write
$$
t_{3}=at_{1}+bt_{2},\qquad t_{4}=ct_{1}+dt_{2}
$$
then these four tritangents define the points $(1:0),(0:1),(a:b),(c:d)$
on $\PP^1$. One of the six cross ratios of these points is
$$
\gamma_{t}:=\frac{
\det\mbox{$\left(
\begin{array}{cc}1&a\\0&b\\
\end{array}\right)$}
\det\mbox{$\left(
\begin{array}{cc}0&c\\1&d
\end{array}\right)$}}
{\det\left(
\begin{array}{cc}1&c\\0&d
\end{array}\right)
\det\left(
\begin{array}{cc}0&a\\1&b
\end{array}\right)}=
\frac{bc}{ad}.
$$

\subsection{Lemma}\label{lemma crosses}
The cross ratio $\gamma_{56}$ defined by
 $t=t_{56}$, $l=b_6$, $t_i=t_{i6}$ and the cross ratio $\gamma_{65}$ defined by
$t=t_{65}$, $l=a_6$, $t_i=t_{6i}$, where $i=1,\ldots,4$, are given by:
$$
\gamma_{56}=\frac{(x_1-x_3)(x_2-x_4)}{(x_1-x_4)(x_2-x_3)},\qquad
\gamma_{65}=\frac{(x_1^2-x_3^2)(x_2^2-x_4^2)}{(x_1^2-x_4^2)(x_2^2-x_3^2)}.
$$

\ts
Let $t_{i6}=\sum_i a_{ij}X_j$ be defined by a determinant as in
section \ref{tritangents}, with the obvious change of indices.
Then $a_{ii}=0$. Thus we must have:
$$
a_{12}a_{21}t_{36}=a_{21}a_{32}t_{16}+a_{12}a_{31}t_{26},\qquad
a_{12}a_{21}t_{46}=a_{21}a_{42}t_{16}+a_{12}a_{41}t_{26},\quad
\text{hence}\quad
\gamma_{56}=\frac{a_{31}a_{42}}{a_{32}a_{41}}.
$$
For $\gamma_{65}$ a similar formula, with $a_{ij}$ replaced by $b_{ij}$,
applies.
\qed

\section{The classifying map}

\subsection{The moduli space of marked cubics}
A point in the moduli space of smooth marked cubic surfaces $\cM^{m}$
can be identified with an isomorphism class $[S,a_1,\ldots,a_6]$ of
a smooth cubic surface with 6 skew lines, the other 21 lines are then
naturally labelled by $b_i,c_{ij}=c_{ji}$, $i,j\in\{1,\ldots,6\}$, $i\neq j$,
cf.\ section \ref{lines}.

The moduli space of smooth marked cubic surfaces $\cM^{m}$ has a
natural compactification $\cM$ (denoted by $\overline{\cM}^m_{\rm cub}$ in
\cite{DGK}, 2.8). The Weyl group $W(E_6)$ acts biregularly on $\cM$
and the quotient $\cM/W(E_6)$ is the GIT quotient of the space of
cubic surfaces in $\PP^3$.

The boundary $\cM-\cM^m$ consists of
 36 (irreducible) boundary divisors which are parametrised by the
positive roots of $E_6$. The divisor corresponding to a root $\alpha$ will
be denoted by $D_\alpha$.
The projective variety $\cM$ is smooth except for
 40 singular points, called the cusps of $\cM$, which form one $W(E_6)$-orbit and
which map to the unique non-stable point in $\cM/W(E_6)$.

A tritangent plane $t=\{l,m,n\}$ with $l,m,n\in\{a_i,b_j,c_{rs}\}$
defines a divisor $D_t$ in $\cM$, the tritangent divisor associated to $t$,
which is the closure of the locus of marked surfaces such that the lines $l,m,n$
meet in a point. Such a point is called an Eckardt point.

\subsection{The line bundle $\cL$ on $\cM$}
The cross ratios associated to a tritangent plane extend to
rational functions on $\cM$. Naruki \cite{Naruki} showed that the rational map
$\cM\rightarrow (\PP^1)^{270}$ defined by the $270=6\cdot 45$ cross ratios
is an embedding on the complement of the 40 cusps and
blows up each cusp to a copy of $(\PP^1)^3$.

Allcock and Freitag \cite{Fr} showed that there is
a very ample line bundle ${\cal L}$ over ${\cal M}$
with the property that any cross ratio is the quotient
of two global sections of $\cL$. The vector space $H^0(\cM,\cL)$ is
$10$-dimensional
and the group $W(E_6)$ acts on $H^0(\cM,\cL)$.
The corresponding representation is the unique irreducible
$10$-dimensional representation of $W(E_6)$, cf.\ \cite{vG}, section 5.

\subsection{Crosses}
The three lines in a tritangent plane $t$ correspond to three weights in
$Q(E_6)\otimes\QQ$, the orthogonal complement in $Q(E_6)$
of these three weights is a root lattice $t^\perp$
of type $D_4$. A cross is a divisor
$$
D_c:=D_t+D_\alpha+D_\beta+D_\gamma+D_\delta
$$
where $\alpha,\ldots,\delta$ are positive, mutually perpendicular,
roots of $E_6$ in $t^\perp$ (cf.\ \cite{Fr}, Definition 3.2).
Each tritangent divisor determines three crosses, so there are $3\cdot 45=135$
crosses. These crosses are linearly equivalent and $\cL\cong\cO_\cM(D_c)$ for
any cross $D_c$. The space $H^0(\cM,\cL)$ is spanned by sections
whose divisors are crosses. Three sections whose zero divisors are the
three crosses associated to a given tritangent $t$ span a
two dimensional subspace of $H^0(\cM,\cL)$ and the quotient of any two
of these three sections is, up to a scalar multiple,
a cross ratio associated to $t$.

\subsection{Divisors and involutions}
The boundary divisor $D_\alpha$ is the fixed point locus of the reflection
$s_\alpha\in W(E_6)$ in $\cM$. The tritangent divisor $D_t$ is the
fixed point locus of the involution $\gamma(t)\in W(E_6)$ which is the the product
of the reflections in (any) four perpendicular roots in $t^\perp$
(cf.\ \cite{Naruki}, Section 8).
In particular, $\gamma(t)=-I$ on $t^\perp$ and $+I$ on the span of the weights
in $t$.

\subsection{The vector space $V$}
Given the marked family of cubic surfaces ${\cal X}\rightarrow U$
(with marking given by the lines $a_1,\ldots,a_6$), we obtain a classifying map
$$
\Phi_\cX:U\longrightarrow \cM.
$$
This map is equivariant for the action of $W(D_5)$.
Let
$$
V=\Phi_\cX^*H^0(\cM,\cL).
$$
The ten dimensional vector space $V$ can be described as follows.

\subsection{Lemma}\label{basis V}
The vector space $V$ is the $10$ dimensional vector space spanned by the following
functions:
$$
\begin{array}{cc}
x_1(x_2^2-x_3^2)(x_4^2-x_5^2),&
x_1(x_2^2-x_5^2)(x_3^2-x_4^2),\\
x_2(x_1^2-x_3^2)(x_4^2-x_5^2),&
x_2(x_1^2-x_4^2)(x_3^2-x_5^2),\\
x_3(x_1^2-x_2^2)(x_4^2-x_5^2),&
x_3(x_1^2-x_4^2)(x_2^2-x_5^2),\\
x_4(x_1^2-x_2^2)(x_3^2-x_5^2),&
x_4(x_1^2-x_3^2)(x_2^2-x_5^2),\\
x_5(x_1^2-x_2^2)(x_3^2-x_4^2),&
x_5(x_1^2-x_4^2)(x_2^2-x_3^2).\\
\end{array}
$$

\ts
The morphism $\Phi_\cX$ extends to a rational map $\PP(Q_\CC)\rightarrow \cM$
which, since $\PP(Q_\CC)$ is smooth, is a morphism
on an open set $U'$ whose complement has codimension at least
two in $\PP(Q_\CC)$. Then $Pic(U')\cong Pic(\PP(Q_\CC))\cong\ZZ$, so the global sections
of $\cL$ pull-back to homogeneous polynomials of a fixed degree.

The map $\Phi_\cX$ extends to the general points
in the hyperplanes $x_i=\pm x_j$, $i\neq j$.
This extension is obviously still $W(D_5)$-equivariant.
The positive roots $x_i- x_j,x_i+x_j$, $1\leq i<j\leq 5$,
are zero on the fixed point locus of the reflections $s_{ij},s_{klm}$
(with $\{i,j,k,l,m\}=\{1,\ldots,5\}$) respectively, so
these 20 hyperplanes are mapped to the corresponding boundary divisors.
As points of $U$ correspond to smooth cubic surfaces, we conclude that
$\Phi_\cX^{-1}(D_\alpha)=(h_\alpha=0)\cap U'$,
where $\alpha\in D_5(\subset E_6)$ is a positive root.

The root system $D_4(\subset D_5\subset E_6)$ with roots $\pm x_i\pm x_j$,
$2\leq i<j\leq 6$, is $t^\perp$ for $t=\{a_6,b_1,c_{16}\}$. The involution
$\gamma_t$ is thus $(x_1,x_2,\ldots,x_5)\mapsto (x_1,-x_2,-x_3,-x_4,-x_5)$
which fixes the divisor $x_1=0$
(and the point $(1:0:0:0:0)$ which is not in $U$).
Thus $\Phi_\cX$ maps $(x_1=0)\cap U$ to $D_t$. We now show that
$\Phi_\cX^{-1}(D_t)=(x_1=0)\cap U$. The lines $b_1$ and $c_{16}$
(cf.\ \ref{27 lines}) meet in the point in $\PP(\cT_x)$
defined by $\langle x,(1,0,0,0,0)\rangle\subset Q_\CC$,
we write simply $b_1\cap c_{16}=(1,0,0,0,0)$. Similarly:
$$
a_6\cap b_1=(x_2x_3x_4x_5,x_1(x_2^3-x_3x_4x_5),x_1(x_3^3-x_2x_4x_5),
x_1(x_4^3-x_2x_3x_5),x_1(x_5^3-x_2x_3x_4)).
$$
It is now easy to verify that $x\in\Phi_\cX^{-1}(D_t)$ iff $x_1=0$.

Thus there is a section of $\cL$, which defines
a cross associated to the tritangent divisor $D_t$, which pulls back to
$x_1^{a_0}(x_2-x_3)^{a_1}(x_2+x_3)^{a_2}(x_4-x_5)^{a_4}(x_4+x_5)^{a_4}$
for certain $a_n\in\ZZ_{>0}$.
Applying $s_{34}\in W(D_4)$, which fixes the tritangent $t$,
we get the pull-back of another section, and the quotient of these two should be a cross ratio
associated to $t$, like $\gamma_{61}$, cf.\ Lemma \ref{lemma crosses}.
So we must have $a_1=a_2=a_3=a_4=1$.

To show that $a_0=1$ if suffices to show that the classifying map is not
ramified at a general point of $x_1=0$.
For this it suffices to exhibit just 4 cross ratios
$\gamma_i$, $i=1,\ldots,4$,  such that the differential
of the map $x\mapsto (\gamma_1(x),\ldots,\gamma_4(x))\in(\PP^1)^4$
has maximal rank in some point $x\in (x_1=0)\cap U$.
We took $\gamma_1=t_{56},\gamma_2=t_{65}$ and $\gamma_3,\gamma_4$
obtained from these two by permuting $x_4\leftrightarrow x_5$,
for the point we took $x=(0:2:3:4:1)$.
We found that the differential is indeed injective at this point.
\qed

\subsection{The $W(D_5)$-action on $V$ and the cross ratios}
\label{functions in V}
Note that the sum of the first two basis functions is:
$$
x_1(x_2^2-x_3^2)(x_4^2-x_5^2)
 +x_1(x_2^2-x_5^2)(x_3^2-x_4^2)=x_1(x_2^2-x_4^2)(x_3^2-x_5^2),
$$
thus all polynomials of the form $x_i(x_j^2-x_k^2)(x_l^2-x_m^2)$,
with $\{i,\ldots,m\}=\{1,\ldots,5\}$,
are contained in $V$.
This verifies that $W(D_5)$ acts on $V$ via its action on the
variables.

Another useful function in $V$ is:
$$
\begin{array}{rcl}
g_{126}&:=&
-\left(x_3(x_1^2-x_2^2)(x_4^2-x_5^2)-x_4(x_1^2-x_2^2)(x_3^2-x_5^2)
 +x_5(x_1^2-x_2^2)(x_3^2-x_4^2)\right)
\\ &=& (x_1^2 - x_2^2)(x_3 - x_4)(x_3 - x_5)(x_4 - x_5).
\end{array}
$$
As $V$ is $W(D_5)$-invariant, it contains all polynomials of the form
$$
(x_i-\epsilon_1 x_j)(x_i-\epsilon_2 x_k)(x_j-\epsilon_3 x_l)
(x_l^2-x_m^2);
\qquad \epsilon_a\in\{1,-1\},\quad
\epsilon_1\epsilon_2\epsilon_3=1,
$$
and
$\{i,j,\ldots,m\}=\{1,2,\ldots,5\}$.

To write a cross like $\gamma_{56}$ as quotient of functions in $V$
we observe first of all that
$$
(x_1^2-x_2^2)(x_3-x_4)(x_3x_4+x_5^2)=
-x_3(x_1^2-x_2^2)(x_4^2-x_5^2)+x_4(x_1^2-x_2^2)(x_3^2-x_5^2)\quad\in V.
$$
Using the $W(D_5)$ action, we then also have $g_{12,34}\in V$ with:
$$
\begin{array}{rcl}
g_{12,34}&:=&
(x_1-x_2)(x_3^2-x_4^2)(x_1x_2+x_5^2) + (x_1^2-x_2^2)(x_3-x_4)(x_3x_4+x_5^2) \\
&=&(x_1-x_2)(x_3-x_4)(x_1x_2x_3 + x_1x_2x_4 + x_1x_3x_4 + x_2x_3x_4
 +(x_1 + x_2 + x_3 + x_4)x_5^2).
\end{array}
$$
Obviously, crosses like $\gamma_{56}$ are quotients of suitable
polynomials in the $W(D_5)$-orbit of $g_{12,34}$.

\section{Extending the classifying map}

\subsection{The rational map}\label{map Phi}
The composition of the classifying map $\Phi_\cX:U\rightarrow \cM$
and the embedding $\phi_\cL:\cM\longrightarrow \PP^9$, defined by the global
section of $\cL$, extends
to a rational map
$$
\Phi:\PP^4=\PP(Q_\CC)\longrightarrow \PP^9.
$$
As coordinate functions of this map one can take
the basis of $V$ given in \ref{basis V}.
We will study the extension
of this map to a suitable blow-up $\tilde{\PP}^4$
of $\PP^4$. The smooth projective variety $\tilde{\PP}^4$
is a compactification of the
open set $U$ from section \ref{definition U}. In particular, this gives a
morphism $\tilde{\PP}^4\rightarrow \cM$ which extends (`compactifies')
the classifying map $\Phi_\cX:U\rightarrow \cM$.

\subsection{The involution $\iota$}\label{iota}
It is obvious that the map $\Phi$ factors over the birational involution:
$$
\iota:\PP^4\longrightarrow \PP^4,\qquad
(x_1:\ldots:x_5)\longmapsto (x_1^{-1}:\ldots:x_5^{-1})
=(x_2x_3x_4x_5:\ldots:x_1x_2x_3x_4)
$$
in fact, the coordinate functions of $\Phi$ are
$x_i(x_j^2-x_k^2)(x_l^2-x_m^2)$ and after putting $x_a:=x_a^{-1}$ and
multiplying all coordinates by $(x_1x_2\ldots x_5)^2$ one obtains again
$x_i(x_j^2-x_k^2)(x_l^2-x_m^2)$.

\subsection{Lemma}\label{fibers Phi}
 The map $\Phi$ has degree two and $\Phi$ is
$\iota$-invariant: $\Phi\circ \iota=\Phi$.

\ts
We already observed that $\Phi\circ \iota=\Phi$ in \ref{iota}.
We study the general fiber of $\Phi$:
$$
\Phi^{-1}(p),\qquad p=(1:a_2,\ldots:a_{10})\in {\rm im}(\Phi).
$$
A point $(x_1:\ldots:x_5)\in\Phi^{-1}(p)$
satisfies the following two equations:
$$
a_5x_1(x_2^2-x_3^2)(x_4^2-x_5^2)=   x_3(x_1^2-x_2^2)(x_4^2-x_5^2),\qquad
a_5x_2(x_1^2-x_3^2)(x_4^2-x_5^2)=a_3x_3(x_1^2-x_2^2)(x_4^2-x_5^2).
$$
Removing the common factor $(x_4^2-x_5^2)$ we obtain that
$(x_1:x_2:x_3)\in\PP^2$ must lie in the intersection of the two cubics:
$$
a_5x_1(x_2^2-x_3^2)-   x_3(x_1^2-x_2^2)=0,\qquad
a_5x_2(x_1^2-x_3^2)-a_3x_3(x_1^2-x_2^2)=0.
$$
There are 7 obvious points of intersection:
$$
(0:0:1),\quad (0:1:0),\quad (1:0:0),\quad (1:1:1),\quad (1:1:-1),\quad
(1:-1:1),\quad (-1:1:1).
$$
The other two points can be found as follows.
Put $x_3=1$.
From the first equation
we have:
$$
x_2^2=x_1(a_5+x_1)/(a_5x_1+1)
$$
substituting this in the second equation gives:
$$
x_2=a_3(x_1^2-x_2^2)/a_5(x_1^2-1)=a_3x_1/(a_5x_1+1).
$$
To have compatibility we need:
$$
x_1(a_5+x_1)/(a_5x_1+1)=\left(a_3x_1/(a_5x_1+1)\right)^2
\quad{\rm hence}\quad
a_5 + (1 - a_3^2 + a_5^2 )x_1 + a_5x_1^2=0
$$
is, essentially, the only remaining equation.
Note that if $x_1$ is a root, then so is $x_1^{-1}$.
We already found $x_2=a_3x_1/(a_5x_1+1)$.
Proceeding in this way with the other equations one finds the result.
\qed

\subsection{Lemma}\label{base locus}
Let $B$ be the base locus of $\Phi$:
$$
B:=\{x\in\PP^4:\; f(x)=0\quad\forall f\in V\,\}.
$$
Then $B$ is the union of 50 lines,
they are the $W(D_5)$-orbits of the lines
$$
l=\{(s:t:t:t:t):\;(s:t)\in\PP^1\,\},\qquad
m=\{(s:t:0:0:0):\;(s:t)\in\PP^1\,\}.
$$
These orbits consist of $40$ and $10$ lines respectively.

\ts
It is easy to check that $l$ and $m$ lie in $B$.
As $B$ is invariant under $W(D_5)$, it suffices to show that each base
point can be mapped into $l\cup m$ using an element of $W(D_5)$.
Let $p=(x_1:\ldots:x_5)$ be a base point.
Assume first that none of the coordinates of $p$ is zero.
Then we must have $x_i=\pm x_j$ for some $i\neq j$ and
using a suitable element of $W(D_5)$ we may assume that $x_4=x_5$.
Next we must still have $x_l=\pm x_m$ with $1\leq l<m\leq 3$, so we
may assume that $x_2=x_3$. As
$x_1(x_2^2-x_4^2)(x_3^2-x_5^2)\in V$
it follows that $x_2=\pm x_4$ or $x_3=\pm x_5$ so using the $W(D_5)$
action again we get $p=(s:t:t:t:t)$ for some $(s:t)\in\PP^1$.

Assume now that at least
one of the coordinates of $p$ is zero.  Using a suitable element
of $W(D_5)$ we may then assume that $x_5=0$ and, using the bases of $V$
given in \ref{basis V} we then have $x_ix_j^2(x_k^2-x_l^2)=0$ whenever
$\{i,j,k,l\}=\{1,2,3,4\}$. If all other coordinates are non-zero, then,
as above, we find that $p$ is in the $W(D_5)$-orbit of $(0:t:t:t:t)\in l$.
If one other coordinate is zero, we may assume that $x_4=x_5=0$ and we
must have $x_ix_j^2x_k^2=0$ whenever $\{i,j,k\}=\{1,2,3\}$.
Thus one more coordinate is zero and we may assume that
$p=(s:t:0:0:0)\in m$.
\qed

\subsection{The singular points of $B$}\label{base points}
The singular points of the base locus $B$ are the points of intersection
of the lines in $B$. There are $21$ such points,
they are the $W(D_5)$-orbits of the points
$$
p_1:=(1:1:1:1:1),\qquad
q_1:=(1:0:0:0:0).
$$
These orbits consist of $16$ and $5$ points respectively.
The points in these orbits are denoted by $p_i$ ($1\leq i\leq 16$)
and $q_i$ ($1\leq i\leq 5$), only when necessary will we define
these other points precisely.

\subsection{The blow-up $\tilde{\PP}^4$}\label{blow up}
In order to extend $\Phi:\PP^4-B\rightarrow \PP^9$ to a morphism,
we first blow-up $\PP^4$ in the $21$ singular points of $B$,
and next we blow-up the strict transforms of the 50 lines in $B$.
The variety which we obtain in this way
is denoted by $\tilde{\PP}^4$.

The strict transform in $\tilde{\PP}^4$ of the exceptional divisor
over $p_i,q_i$ is denoted by $E_{p_i}$ and $E_{q_i}$ respectively,
they are $\PP^3$'s blown up in a finite number of points, the points where
the strict transforms of the lines in $B$ meet the exceptional divisors
in the first blow-up.

The lines in the $W(D_5)$-orbit of $m,l$ will be denoted by
$m_\alpha,l_\beta$, $1\leq\alpha\leq 40$, $1\leq\beta\leq 10$.
The exceptional divisors over the lines $m_\alpha$ and $l_\beta$
are denoted by $E_{m_\alpha}$ and $E_{l_\beta}$.

Explicit calculations,
which are easy and which we omit (but see the next section), verify that
$\Phi$ extends to a morphism, which we denote by the same name,
$$
\Phi:\tilde{\PP}^4\longrightarrow \PP^9.
$$

\section{Divisors}
\subsection{Contractions}\label{contractions}
The divisors $E_{l_\alpha}$ and $E_{m_\beta}$, which are birational
to $\PP^2$-bundles over the (strict transforms of the)
lines $l_\alpha$ and $m_\beta$ respectively,
are contracted to planes in $\PP^9$. In fact,
the image of a point $(x:y:z)$ in the fiber over
$(1:a:a:a:a)\in l$, for $a$ general,
is given by:
$$
\lim_{t\rightarrow 0}\Phi(1:a:a+tx:a+ty:a+tz)=
(0:0:y-z:x-z:y-z:-z:x-z:-z:x-y:-x),
$$
which does not depend on $a$, so $\Phi$ collapses all fibers of the blow-up
$E_{l}\rightarrow l$ to a fixed $\PP^2\subset\PP^9$.
Similarly, the image of a point $(x:y:z)$ in the fiber over $(1:a:0:0:0)\in m$
is given by:
$$
\lim_{t\rightarrow 0}\Phi(1:a:tx:ty:tz)=(0:0:0:0:0:x:0:y:0:z).
$$

\subsection{Tritangent divisors}
In the proof of Lemma \ref{basis V} we observed that the classifying map
$\Phi_\cX$
maps the $\PP^3\subset \PP^4$
defined by $x_i=0$ ($1\leq i\leq 5$) onto a tritangent divisor.
The image under $\Phi$ of such a $\PP^3$ spans a $\PP^7\subset\PP^9$.
The covering involution
maps this $\PP^3$ to the point $q_i$ which has all coordinates
equal to zero except $x_i=1$.
Let $i=5$, then the map $\Phi$ induces the following map on the
the exceptional divisor $E_{q_5}$, which is birational to $\PP^3$, over $q_5$:
$$
\Phi:E_{q_5}\sim\PP^3\longrightarrow \PP^7\subset\PP^9,\quad
(y_1:y_2:y_3:y_4)\longmapsto \lim_{t\rightarrow0}\Phi(ty_1:ty_2:ty_3:ty_4:1).
$$
It is easy to verify that the image of $(y_1:\ldots:y_4)$ is given by:
$$
(y_1(y_2^2 - y_3^2): y_1(y_3^2 - y_4^2):
\ldots: y_4(y_1^2 - y_2^2): y_4(y_1^2 - y_3^2): 0: 0).
$$
Thus on this $\PP^3$ the map $\Phi$ coincides with the map
on $\PP^3_{\rm w}$ defined in terms of the root system $F_4$
considered in \cite{vG}, Thm.\ 6.5.

In particular, the inverse image of the tritangent divisor
labelled by $\{a_6,b_i,c_{i6}\}$, with $1\leq i\leq 5$,
has two irreducible components in $\tilde{\PP}^4$, the strict transforms
of the hyperplane $x_i=0$ and the divisor $E_{q_i}$,
which are exchanged by $\iota$.

The remaining 40  tritangent divisors are components of the divisors of the
other crosses, so they are
the images of the 40 cubics in
the $W(D_5)$-orbit of (cf.\ \ref{functions in V}):
$$
X:\quad x_1x_2x_3 + x_1x_2x_4 + x_1x_3x_4 + x_2x_3x_4
 +(x_1 + x_2 + x_3 + x_4)x_5^2=
$$
$$
=x_1x_2x_3x_4(x_1^{-1}+\ldots+x_4^{-1})+
(x_1 + x_2 + x_3 + x_4)x_5^2=0.
$$
From the equation it is obvious that $X$ is invariant under $\iota$.

Note that $X$ (and its $W(D_5)$-conjugates) are birationally covers of a
$\PP^3$ (with coordinates $(x_1:x_2:x_3:x_4)$) branched over
the union of a plane $x_1+x_2+x_3+x_4=0$ and the 4-nodal
Cayley cubic surface
$x_1x_2x_3 + x_1x_2x_4 + x_1x_3x_4 + x_2x_3x_4=0$.
In particular, $X$ is singular and projection from a singular point of $X$
gives a birational isomorphism between $X$ and $\PP^3$.

\subsection{Boundary divisors}\label{boundary}
There is a $W(E_6)$-equivariant bijection between
the $36$ boundary divisors of $\cM$ and the $36$ positive roots
in $E_6$. There are two $W(D_5)$-orbits. One consists of the 20
roots of $h_{ij},h_{ijk}$, $i,j,k\leq 5$, which are perpendicular to
$a_6$ and these are the positive roots of $D_5$.
Note that the involution $\iota$ fixes each of the divisors $h_{ij}=0$,
$h_{ijk}=0$, hence $\Phi$ is a degree two map on the strict transform
of such a divisor.

The other orbit consists
of the 16 roots $h,h_{i6},h_{ij6}$, $i,j\leq 5$.
We already determined the image of each of the divisors in $\tilde{\PP}^4-U$
except for the divisors $E_{p_i}$, $1\leq i\leq 16$. Thus these must map
to the remaining 16 boundary divisors.
On the divisor $E_{p_1}$,
which is birational to the exceptional divisor over $p_1$ in the blow-up of
$\PP^4$ in $p_1$,
the map induced by $\Phi$ is
$$
\Phi:E_{p_1}\cong\PP^3\longrightarrow \PP^4\subset\PP^9,\quad
(y_1:y_2:y_3:y_4)\longmapsto \lim_{t\rightarrow0}
\Phi(ty_1+1:ty_2+1:ty_3+1:ty_4+1:1)
$$
and it is easy to see that the image of $(y_1:\ldots:y_4)$ is given by:
$$
(y_4(y_2 - y_3): y_2(y_3 - y_4):
\ldots : (y_1 - y_2)(y_3 - y_4):(y_2 - y_3)(y_1 - y_4)).
$$
The 10 coordinate functions span the space of quadrics which vanish in
the 5 points $(1:0:0:0),\ldots,(0:0:0:1),(1:1:1:1)$.
It is well-known that this map gives a birational
isomorphism of $\PP^3$
with the Segre threefold in $\PP^4$, cf.\ \cite{Hunt}, Theorem 3.2.1.
It is well-known that a boundary divisor of $\cM$ is isomorphic to a
Segre cubic threefold. The nodes of the Segre threefold are cusps of $\cM$.
The line parametrized by $(0:0:s:t)$ maps to a node of the cubic
and explicit computation shows that this cusp
has coordinates $(1:0:1:1:0:0:0:0:0:-1)\in\cM\subset\PP^9$.

The inverse image of any of these 16 boundary divisors must be an $E_{p_i}$
since we accounted for all other divisors in $\tilde{\PP}^4$. As $\Phi$ is
a birational morphism on a $E_{p_i}$, these divisors lie in the ramification
locus of $\Phi$.

\subsection{Cusps}\label{cusps}
The inverse image of a cusp of $\cM$ under $\Phi$ is a plane in $\PP^4$.
For example, the inverse image of the cusp $(1:0:1:1:0:0:0:0:0:-1)$
is easily seen to be the plane $\gamma$ parametrized by
$$
\gamma:=\{(s:s:t:u:s):\;(u:s:t)\in\PP^2\,\},\qquad\Phi(\gamma)=
(1:0:1:1:0:0:0:0:0:-1).
$$
There are 40 planes in the $W(D_5)$-orbit of $\gamma$, these map to the
40 cusps of $\cM\subset \PP^9$.

\subsection{Labelling of divisors}
Naruki's cross ratio variety $\cC$ is a natural
desingularisation of the moduli space $\cM$. The map $\cC\rightarrow \cM$
contracts 40 disjoint divisors in $\cC$,
each isomorphic to $(\PP^1)^3$, to the 40 cusps of $\cM$ and is an isomorphism
on the complements. We refer to these 40 divisors on $\cC$ as cusp divisors.

The three types of divisors on $\cC$
are conveniently labelled by points in a $4$-dimensional
projective space over $\FF_3$, canonically, it is $\PP(P(E_6)/Q(E_6))$,
cf.\ \cite{CvG}, 1.7, 7.4-7.6. It is remarkable that these labels allow one
to identify the corresponding loci in $\PP^4$ easily.
We will content ourselves with some examples.

The boundary divisors are parametrized
by the 36 points obtained from $(0:0:0:1:1)$ and $(1:1:1:1:1)$
by permutation and sign changes of the coordinates. The boundary
divisor labelled by $(0:0:0:1:-1)$ is the image of the
$\PP^3\subset \PP(Q_\CC)$
defined by $h_{45}=x_4-x_5=0$.
The boundary divisor labelled by $(1:-1:1:-1:1)$
is the image of the exceptional
divisor over the singular point $(1:-1:1:-1:1)$ in the base locus of $\Phi$.

The tritangent divisors are parametrized by the 45 points obtained from
$(0:0:0:0:1)$ and $(0:1:1:1:1)$ by permutation and
sign changes of the coordinates. The point $(0:0:0:0:1)$ labels the boundary
divisor which is the image of the exceptional divisor over
$(0:0:0:0:1)\in\PP(Q_\CC)$, equivalently, it is the image of the hyperplane
$x_5=0$. The point $(0:1:-1:1:1)$ labels the tritangent divisor
which is the image of the cubic
$x_2x_3x_4x_5(x_2^{-1}-x_3^{-1}+x_4^{-1}+x_5^{-1})+
(x_2 - x_3 + x_4 + x_5)x_1^2=0$.

The cusp divisors are parametrized by the
40 points in $\PP^4(\FF_3)$ obtained from $(0:0:1:1:1)$ by permutation
and sign changes of the coordinates. The plane in $\PP(Q_\CC)$
which maps to the cusp labelled by $(0:0:1:-1:1)$
is parametrized by $(s:t:u:-u:u)$.
The plane $\gamma\subset \PP(Q_\CC)$ in \ref{cusps}
above maps to the cusp labelled $(1:1:0:0:1)$.

\section{The Chow ring}
\subsection{} In the paper \cite{CvG} we determined the Chow ring of
Naruki's cross ratio variety $\cC$ which is a desingularisation of the
moduli space $\cM$ of marked cubic surfaces. We will now see that the
extension of the classifying map, which we denote again by
$\Phi_{\cX}:\tilde\PP^4\longrightarrow \cM$ (cf.\ sections \ref{map Phi},
 \ref{blow up})
allows us to recover quickly some of these results.

The singular locus of $\cM$ consists of the 40 cusps which are
in one $W(D_5)$-orbit.
The inverse image
of one cusp is the surface $\tilde{\gamma}\subset\tilde\PP^4$
which is the strict transform of the
plane $\gamma\subset\PP^4$, cf.\ \ref{cusps}.
We will thus be interested in the open subset
$$
(\tilde{\PP}^4)_0:=\tilde{\PP}^4-\cup\tilde{\delta},\qquad{\rm and}\quad
\Phi_0:=(\Phi_\cX)_{|(\tilde{\PP}^4)_0}:
       (\tilde{\PP}^4)_0\longrightarrow \cM^{sm},
$$
where the union is over the 40 surfaces in the $W(D_5)$-orbit of
$\tilde{\gamma}$.
The variety $(\tilde{\PP}^4)_0$  maps, generically 2:1,
onto the smooth locus $\cM^{sm}$ of $\cM$. The smooth, quasi-projective
variety $\cM^{sm}$ is the moduli space of marked nodal cubic surfaces
(\cite{DGK}, section 2.8).

\subsection{Proposition}\label{A1}
The Picard group, tensored by $\QQ$,
of the moduli space
of marked nodal cubic surfaces $\cM^{sm}$ is:
$$
A^1(\cM^{sm})_\QQ
\cong \left(\oplus_{i=1}^{16} \QQ B_{\alpha_i}\right) \,\oplus\,
\left(\oplus_{j=1}^5 \QQ T_j\right)
$$
where
$$
B_{\alpha_i}:=(\Phi_{0})_*E_{p_i},\qquad
T_j:=(\Phi_{0})_*E_{q_j}.
$$
Here $B_{\alpha_i}$ is the boundary divisor parametrised by
$\alpha_i$ which runs over the $16$ positive roots of $E_6$
which are {\em not} perpendicular to the weight $a_6$
(cf.\ section \ref{boundary}) and $T_j$ is the tritangent divisor
labelled by $\{a_6,b_j,c_{j6}\}$.

\ts
Basic results on the Chow group of a blow-up imply:
$$
A^1((\tilde{\PP}^4)_0)=
\ZZ H\oplus (\oplus^{16} E_{p_i})\oplus (\oplus^5 E_{q_i})
\oplus (\oplus^{10} E_{m_\alpha})\oplus (\oplus^{40} E_{l_\beta})
$$
where $H$ is the class of the strict transform of a hyperplane in $\PP^4$.
The divisors $E_{l_\beta},E_{m_\alpha}$ are contracted by $\Phi_0$
(cf.\ section \ref{contractions}), hence
$(\Phi_0)_*E_{l_\beta}=(\Phi_0)_*E_{m_\alpha}=0$.

The strict transform $D_i\subset (\tilde{\PP^4})_0$
of the hyperplane $x_i=0$ in $\PP^4$
is mapped birationally
onto a tritangent divisor $T_i$, so $(\Phi_0)_*D_i=T_i$.
The same tritangent divisor is also the birational
image of $E_{q_i}$, $(\Phi_0)_*E_{q_i}=T_i$.
As there are 4 singular points, $q_j$ with $j\neq i$,
of the base locus in $x_i=0$, we have:
$$
D_i=H-\sum_{j\neq i} E_{q_j},\qquad{\rm thus}\quad
(\Phi_0)_*H=(\Phi_0)_*(D_i+\sum_{j\neq i} E_{q_j})=\sum_{j=1}^5 T_j.
$$
Hence
$(\Phi_0)_*A^1((\tilde{\PP}^4)_0)\cong
(\oplus^{16}B_{\alpha_i})\oplus(\oplus^5T_j)$.
As $(\Phi_0)_*(\Phi_0)^*$ is multiplication by $2$ on $A^1(\cM^{sm})$,
the proposition follows.
\qed

\subsection{} The group $W(E_6)$ acts biregularly on $\cM^{sm}$ and $\cC$.
It permutes the $45$ tritangent and the $36$ boundary divisors
transitively.
The sum of the tritangent divisors and of the
boundary divisors are denoted by
$\hat{T}$ and $\hat{B}$ respectively. These $W(E_6)$ invariant classes
on $\cM^{sm}$ are related as follows.

\subsection{Proposition}\label{tb}
In the Picard group of $\cM^{sm}$ we have:
$$
4\hat{T}=25\hat{B}.
$$

\ts
The strict transform $\tilde{H}_{12}$ in $\tilde{\PP}^4$
of the hyperplane $x_1-x_2=0$ in $\PP^4$ is mapped 2:1
onto a boundary divisor $B_{h_{12}}$
in $\cM^{sm}$ by $\Phi_0$.
As this hyperplane contains exactly $8$ points of the base
locus $B$ whose exceptional divisor maps birationally
onto a boundary divisor and
exactly 3 points, $q_3,q_4,q_5$ whose exceptional divisor maps to a
tritangent divisor, we get
$\tilde{H}_{12}=H-\sum^8 E_{p_i}-\sum^3 E_{q_i}$+\ldots,
where we omit the precise indices
and divisors which are contracted by $\Phi_0$.
In the proof of Prop.\ \ref{A1} we saw that
$(\Phi_0)_*H=\sum_{j=1}^5 T_j$.
As $(\Phi_0)_*\tilde{H}_{12}=2B_{h_{12}}$,
we find that a sum of ten boundary divisors
is linearly equivalent to a sum of two tritangent divisors in $\cM^{sm}$.
If $B$ is any boundary divisor,
then $\sum_\sigma \sigma^*B=(|W(E_6)|/36)\hat{B}$, since
$W(E_6)$ is transitive on the set of $36$ boundary divisors.
Similarly $\sum_\sigma \sigma^*T=(|W(E_6)|/45)\hat{T}$
for any tritangent divisor $T$.
In this way the linear equivalence of the sum of ten boundary divisors with
the sum of two tritangent
divisors leads to
$\mbox{$\frac{10}{36}$}\hat{B}=\mbox{$\frac{2}{45}$}\hat{T}$.
\qed

\subsection{Remark} In \cite{CvG}, Theorem 2, we proved that Naruki's smooth
compactification $\cC$ of $\cM^{sm}$ has $A^1(\cC)\cong \ZZ^{61}$.
As $\cC-\cM^{sm}$ is the disjoint union of $40$ cusp divisors,
whose classes in $\cC$ generate a direct summand isomorphic to
$\ZZ^{40}$, we see that Proposition \ref{A1}
is consistent with the results of \cite{CvG}.

Proposition \ref{tb} follows from the relation (\cite{CvG}, Theorem 2.4(2)):
$$
\hat{T}=(25\hat{B}+27\hat{C})/4\qquad(\in A^1({\cal C})),
$$
because $\hat{\cal C}$, the sum of the 40 cusp divisors,
is zero on $\cM^{sm}$.

\


\end{document}